\documentclass{article}
\usepackage{amssymb}
\usepackage{graphicx}
\newtheorem{theorem}{Theorem}

\newtheorem{example}{Example}
\newtheorem{definition}{Definition}

\newenvironment{proof}[1][Proof]{\textbf{#1.} }{\ \rule{0.5em}{0.5em}}

\title{Algebraic operations on the space of Hausdorff continuous
interval functions}

\author{Roumen Anguelov\\
Department of Mathematics and Applied Mathematics\\
University of Pretoria, SOUTH AFRICA\\
roumen.anguelov@up.ac.za\\[3pt]
Blagovest Sendov\\
Bulgarian Embassy, 36-3,Yoyogi 5-chome,\\ Shibuia-ku, Tokyo
151-0053 JAPAN\\
sendov2003@yahoo.com\\[3pt]
Svetoslav Markov\\
Institute of Mathematics and Informatics, BAS \\
"Acad. G. Bonchev" st., block 8, 1113 Sofia, BULGARIA\\
smarkov@bio.bas.bg}
\date{}

\begin{document}
\maketitle

\begin{abstract}
We show that the operations addition and multiplication on the set
$C(\Omega)$ of all real continuous functions on
$\Omega\subseteq\mathbb{R}^n$ can be extended to the set
$\mathbb{H}(\Omega)$ of all Hausdorff continuous interval
functions on $\Omega$ in such a way that the algebraic structure
of $C(\Omega)$ is preserved, namely, $\mathbb{H}(\Omega)$ is a
commutative ring with identity. The operations on
$\mathbb{H}(\Omega)$ are defined in three different but equivalent
ways. This allow us to look at these operations from different
points of view as well as to show that they are naturally
associated with the Hausdorff continuous functions.
\end{abstract}

\section{Introduction} The set $\mathbb{H}(\Omega)$ of all Hausdorff continuous interval
functions appears naturally in the context of Hausdorff
approximations of real functions. The concept of Hausdorff
continuity of interval functions generalizes the concept of
continuity of real functio  n in such a way that many essential
properties of the usual continuous real functions are preserved.
Not least this is due to the fact that the Hausdorff continuous
functions assume real (point) values on a dense subset of the
domain and are completely determined by these values. It is well
known that the set $C(\Omega)$ of all continuous real functions
defined on a subset $\Omega$ of $\mathbb{R}^n$ is a commutative
ring with respect to the point-wise defined addition and
multiplication of functions. Hence the natural question: Is it
possible to extend the algebraic operations on $C(\Omega)$ to the
set $\mathbb{H}(\Omega)$ of all Hausdorff continuous functions
defined on $\Omega$ in a way which preserves the algebraic
structure, that is, the set of $\mathbb{H}(\Omega)$ is a
commutative ring with respect to the extended operations? We show
in this paper that the answer to this question is positive.
Furthermore, we give three different but equivalent ways of
defining algebraic operations on $\mathbb{H}(\Omega)$ with the
required properties. Namely, (i) through the point-wise interval
operations, (ii) by using an extension property of the Hausdorff
continuous functions, (iii) through the order convergence
structure on $\mathbb{H}(\Omega)$.

\section{General Setting}

The real line is denoted by $\mathbb{R}$ and the set of all finite
real intervals by $\mathbb{IR}=\{[\underline{a},\overline
{a}]:\underline{a},\overline{a}\in\mathbb{R}$,
 $\underline{a}\leq\overline{a}\}$. Given an interval
 $a=[\underline{a},\overline{a}]= \{x :
\underline{a}\leq x\leq\overline{a}\}  \in\mathbb{IR}$,
 $w(a)=\overline{a}-\underline{a}$
is the width of $a$, while
$|a|=\max\{|\underline{a}|,|\overline{a}|\}$ is the modulus of
$a$. An interval $a$ is called proper interval, if $w(a)>0$ and
point interval, if $w(a)=0$. Identifying $a\in \mathbb{R}$ with
the point interval $[a,a]\in \mathbb{IR}$, we consider
$\mathbb{R}$ as a subset of $\mathbb{IR}$. We denote by
$\mathbb{A}(\Omega)$ the set of all locally bounded
interval-valued functions defined on an arbitrary set $\Omega
\subseteq \mathbb{R}^{n}$. The set $\mathbb{A}(\Omega )$ contains
the set $\mathcal{A}(\Omega)$ of all locally bounded real
functions defined on $\Omega$. Recall that a real function or an
interval-valued function $f$ defined on $\Omega$ is called locally
bounded if for every $x\in\Omega$ there exist $\delta>0$ and
$M\in\mathbb{R}$ such that $|f(y)|< M,\ y\in B_\delta(x)$, where
$B_{\delta }(x)=\{y\in \Omega :||x-y||<\delta \}$ denotes the open
$\delta$-neighborhood of $x$ in $\Omega $.

Let $D$ be a dense subset of $\Omega$.  The mappings
$I(D,\Omega,\cdot),$ $S(D,\Omega,\cdot): \mathbb{A}(D
)\longrightarrow \mathcal{A}(\Omega )$ defined for $f \in
\mathbb{A}(D)$ and $x\in \Omega$  by
\begin{eqnarray*}
I(D,\Omega ,f)(x) &=&\sup_{\delta >0}\inf \{ {f}(y):y\in
B_{\delta }(x)\cap D\}, \\
S(D,\Omega ,f)(x) &=&\inf_{\delta
>0}\sup \{ {f}(y):y\in
B_{\delta }(x)\cap D\},
\end{eqnarray*}
are called lower and upper Baire operators, respectively. The
mapping $F: \mathbb{A}(D)\longrightarrow \mathbb{A}(\Omega )$,
called {\em graph completion operator}, is defined by
\[
F(D,\Omega ,f)(x)=[I(D,\Omega ,f)(x),S(D,\Omega ,f)(x)],\;x\in
\Omega ,\;f\in \mathbb{A}(D ). \label{gcgen}
\]
In the case when $D=\Omega$ the sets $D$ and $\Omega$ will be
omitted, thus we write $I(f)=I(\Omega,\Omega,f),\
S(f)=S(\Omega,\Omega,f),\ F(f)=F(\Omega,\Omega,f)$.

\begin{definition}
A function $f\in\mathbb{A}(\Omega)$ is S-continuous, if $F(f)=f$.
\end{definition}

\begin{definition}\label{DefHcont}
A function $f\in \mathbb{A}(\Omega)$ is Hausdorff continuous
(H-continuous), if $g\in\mathbb{A}(\Omega)$ with $g(x)\subseteq
f(x),$ $x\in \Omega$, implies $F(g)(x)=f(x),$ $x\in \Omega $.
\end{definition}

\begin{theorem}\label{tFISf} \cite{Anguelov}, \cite{Sendov}
For every $f\in\mathbb{A}(\Omega)$ the functions $F(I(S(f)))$ and
$F(S(I(f)))$ are H-continuous.
\end{theorem}

The following theorem states an essential property of the
continuous functions which is preserved by the H-continuity,
\cite{Anguelov}.

\begin{theorem}\label{tident}
Let $f,g\in\mathbb{H}(\Omega)$. If there exists a dense subset $D$
of $\Omega$ such that $f(x)=g(x)$, $x\in D$, then $f(x)=g(x)$,
$x\in\omega$.
\end{theorem}

H-continuous functions are also similar to the usual continuous
real functions in that they assume point values everywhere on
$\Omega$ except for a set of first Baire category. More precisely,
it is shown in \cite{Anguelov} that for every
$f\in\mathbb{H}(\Omega)$ the set
\begin{equation}\label{defWf}
W_f=\{x\in\Omega : w(f(x))>0\}
\end{equation}
is of first Baire category and $f$ is continuous on
$\Omega\setminus W_f$. Since a finite or countable union of sets
of first Baire category is also a set of first Baire category we
have:
\begin{theorem}\label{tDfdense}
Let the set $\Omega$ be open and let ${\cal F}$ be a finite or
countable set of H-continuous functions. Then the set
\begin{equation}\label{DF}
D_{\cal F}=\{x\in\Omega : w(f(x))=0,\;f\in{\cal
F}\}=\Omega\setminus \displaystyle\bigcup_{f\in{\cal F}}W_f
\end{equation}
is dense in $\Omega$ and all functions $f\in{\cal F}$ are
continuous on $D_{\cal F}$.
\end{theorem}

The graph completion operator is inclusion isotone i) w. r. t. the
functional argument, that is, if $f,g\in \mathbb{A}(D)$, where $D$
is dense in $\Omega$, then
\begin{equation}\label{Finclmon1}
f(x)\subseteq g(x),\;x\in D\; \Longrightarrow\; F(D,\Omega
,f)(x)\subseteq F(D,\Omega ,g)(x),\;x\in \Omega ,
\end{equation}
and, ii) w. r. t. the set $D$ in the sense that if $D_{1}$ and
$D_{2}$ are dense subsets of $\Omega $ and $f\in
\mathbb{A}(D_{1}\cup D_{2})$ then
\begin{equation}\label{Finclmon2}
D_{1}\subseteq D_{2}\; \Longrightarrow \; F(D_{1},\Omega
,f)(x)\subseteq F(D_{2},\Omega ,f)(x), \; x\in \Omega .
\end{equation}
In particular, (\ref{Finclmon2}) implies that for any dense subset
$D$ of $\Omega $ and $f\in \mathbb{A}(\Omega )$ we have
\begin{equation}\label{Finclmon3}
F(D,\Omega ,f)(x)\subseteq F(f)(x), \; x\in \Omega .
\end{equation}

Let $f\in \mathbb{A}(\Omega).$ For every $x\in \Omega$ the value
of $f$ is an interval $[\underline{f}(x),\overline{f}(x)]\in
\mathbb{IR}$. Hence,  $f$ can be written in the form
$f=[\underline{f},\overline{f}]$ where $\underline{f},\overline{f
}\in\mathcal{A}(\Omega)$ and $\underline{f}(x)\leq
\overline{f}(x),\;x\in\Omega$. The lower and upper Baire operators
as well as the graph completion operator of an interval-valued
function $f$ can be represented in terms of $\underline{f}$ and
$\overline{f}$, namely, for every dense subset $D$ of $\Omega $:
$I(D,\Omega,f)=I(D,\Omega,\underline{f})$,
$S(D,\Omega,f)=S(D,\Omega,\overline{f})$,
$F(D,\Omega,f)=[I(D,\Omega,\underline{f}),S(D,\Omega,\overline{f})]$.

\section{Interval operations and operations on $\mathbb{H}(\Omega)$}

We recall that the operations of addition and multiplication on
the set of real intervals $\mathbb{IR}$ are defined in more then
one way, \cite{??Markov}. Here we consider the so called outer
operations which are inclusion isotone. For
$[\underline{a},\overline{a}],
[\underline{b},\overline{b}]\in\mathbb{I\,}{\mathbb{R}}$
\begin{eqnarray*}
&&\hspace{-7mm}[\underline{a},\overline{a}] +
[\underline{b},\overline{b}]\! =\! \{a+b:a\!\in\!
[\underline{a},\overline{a}], b\!\in\!
[\underline{b},\overline{b}]\}\!=\![\underline{a} +
\underline{b},\overline{a} + \overline{b}]
\\[0pt]
&&\hspace{-7mm}[\underline{a}, \overline{a}]\!\times\!
[\underline{b},\overline{b}]\!=\! \{ab\!:\!a\!\in\!
[\underline{a},\overline{a}], b\!\in\!
[\underline{b},\overline{b}]\}\!=\!
[\min\{\underline{a}\underline{b},\underline{a}\overline{b},\overline{a}\underline{b},\overline{a}\overline{b}\},
\max\{\underline{a}\underline{b},\underline{a}\overline{b},\overline{a}\underline{b},\overline{a}\overline{b}\}]
\end{eqnarray*}
The operations for interval functions are defined point-wise in
the usual way:
\begin{equation}\label{f+g,fxg}
(f+g)(x)=f(x)+g(x)\ ,\ \  (f\times g)(x)=f(x)\times g(x)
\end{equation}

\begin{example}\label{ExNotHcont}
Consider the functions $f,g\in \mathbb{H}(\mathbb{R})$ given by
\[
f(x)=\left\{
\begin{tabular}{cll}
$0$, & if & $x<0$, \\ [4pt] $\lbrack 0,1]$, & if & $x=0$, \\ [4pt]
$1$, & if & $x>0$;
\end{tabular}
\right. \ \ \ \ g(x)=\left\{
\begin{tabular}{cll}
$0$, & if & $x<0$, \\ [4pt] $\lbrack -1,0]$, & if &
$x=0$, \\
[4pt] $-1$, & if & $x>0$.
\end{tabular}
\right.
\]
   Using (\ref{f+g,fxg}) we have
\begin{eqnarray*}
(f+g)(x)=\left\{
\begin{tabular}{cll}
$0$, & if & $x<0$, \\ [4pt] $\lbrack -1,1]$, & if &
$x=0$, \\
[4pt] $0$, & if & $x>0$.
\end{tabular}
\right.
\end{eqnarray*}
\end{example}

The considered example shows that the point-wise sum of
H-continuous functions is not necessarily an H-continuous
function. Hence the significance of the following two theorems.

\begin{theorem}\label{tf*gScont}\cite{AnguelovMarkov}
If the interval functions $f,g\in\mathbb{A}(\Omega)$ are
S-continuous then the function $f+g$ and $f\times g$ are also
S-continuous.
\end{theorem}

\begin{theorem}\label{texistsunique}
Let $f,g\in\mathbb{H}(\Omega)$.
\begin{itemize}
\item[(a)] There exists a unique function $p\in\mathbb{H}(\Omega)$
satisfying the inclusion \linebreak$p(x)\subseteq (f+g)(x)$,
$x\in\Omega$. \item[(b)] There exists a unique function
$q\in\mathbb{H}(\Omega)$ satisfying the inclusion \linebreak
$q(x)\subseteq (f\times g)(x)$, $x\in\Omega$.
\end{itemize}
\end{theorem}
\begin{proof}
We will prove only (a) because (b) is proved in a similar way. The
existence of the function $p$ follows from Theorem \ref{tFISf}.
Indeed, both functions $F(I(S(f+g)))$ and $F(S(I(f+g)))$ satisfy
the required inclusion. To prove the uniqueness let us assume that
$p_1,p_2\in\mathbb{H}(\Omega)$ both satisfy the inclusion in (a).
Consider the set $D_{fg}=\{x\in\Omega:w(f(x))=w(g(x))=0\}$.
Obviously $w((f+g)(x))=0$ for $x\in D_{fg}$. Therefore, due to the
assumed inclusions we have $p_1(x)=p_2(x)=f(x)+g(x)$, $x\in
D_{fg}$. According to Theorem \ref{tFISf} the set $D_{fg}$ is
dense in $\Omega$. Using that the functions $p_1$ and $p_2$ are
H-continuous it follows from Theorem \ref{tident} that $p_1=p_2$.
\end{proof}

Now we define the operations $\oplus$ and $\otimes$ as follows.

\begin{definition}\label{defoper1}
Let $f,g\in\mathbb{H})(\Omega)$. Then $f\oplus g$ is the unique
Hausdorff continuous function satisfying the inclusion $(f\oplus
g)(x)\subseteq (f+g)(x)$, $x\in\Omega$; $f\otimes g$ is the unique
Hausdorff continuous function satisfying the inclusion $(f\otimes
g)(x)\subseteq (f\times g)(x)$, $x\in\Omega$.
\end{definition}

The existence of both $f\oplus g$ and $f\otimes g$ is guaranteed
by Theorem \ref{texistsunique}.

It is important to note that the values of $f\oplus g$ and
$f\otimes g$ at the points where both operands assume interval
values can not be determined point-wise, i.e., from the values of
$f$ and $g$ at these points. This is illustrated by the following
example.

\begin{example}\label{ExNotPointwise}
Consider the functions $f,g\in \mathbb{H}(\mathbb{R})$
given by%
\begin{eqnarray*}
f(x)&=&\left\{
\begin{tabular}{lll}
$\displaystyle\sin (1/x)$, & if & $x\neq 0$, \\
[6pt] $\lbrack -1,1]$, & if & $x=0$;
\end{tabular}
\right.
\\
g(x)&=&\left\{
\begin{tabular}{lll}
$\displaystyle\cos (1/x)$, & if & $x\neq 0$, \\
[6pt] $\lbrack -1,1]$, & if & $x=0$.
\end{tabular}
\right.
\end{eqnarray*}
We have
\[
(f\oplus g)(x)=\left\{
\begin{tabular}{lll}
$\displaystyle\sqrt{2}\cos (1/x + \pi /4 ) $, & if &
$x\neq 0$, \\
[6pt] $\lbrack -\sqrt{2},\sqrt{2}]$, & if & $x=0$.
\end{tabular}%
\right.
\]
Clearly $(f\oplus g)(0)$ can not be obtained just from the values
$f(0)$ and $g(0)$.
\end{example}

According to Theorem \ref{tFISf} for any
$f,g\in\mathbb{H})(\Omega)$ the functions $F(S(I(f+g)))$ and
$F(I(S(f+g)))$ are Hausdorff continuous. Moreover, these functions
satisfy the inclusions
\[F(S(I(f+g)))(x)\subseteq (f+g)(x),\ F(I(S(f+g))(x)\subseteq
(f+g)(x),\ x\in\Omega.
\]
Therefore they both coincide with $f\oplus g$. Hence we have
$f\oplus g=F(I(S(f+g)))=F(I(S(f+g))$. In a similar way we obtain
$f\otimes g=F(I(S(f\times g)))=F(I(S(f\times g))$. One can
immediately see from the above representations that if $f$ and $g$
are usual continuous real functions we have $f\oplus g=f+g$ and
$f\otimes g=f\times g$. Hence $\oplus$ and $\otimes$ extend the
operations of addition and multiplication on $C(\Omega)$. A
further motivation for considering these operations is in the fact
that the algebraic structure of $C(\Omega)$ is preserved as stated
in the next theorem.

\begin{theorem}\label{tring}
The set $\mathbb{H}(\Omega)$ is a commutative ring with identity
with respect to the operations $\oplus$ and $\otimes$.
\end{theorem}
\begin{proof}
The commutative laws for both $\oplus$ and $\otimes$ follow
immediately from Definition \ref{defoper1}. It is also obvious
that the additive identity is the constant zero function while the
multiplicative identity is the constant function equal to 1. We
will show now the existence of the additive inverse. Let
$f=[\underline{f},\overline{f}]\in\mathbb{H}(\Omega)$. Consider
the function $g\in\mathbb{H}(\Omega)$ given by
$g(x)=[-\overline{f}(x),-\underline{f}(x)]$, $x\in\Omega$. Clearly
we have
\[
0\in f(x)+g(x),\ x\in\Omega.
\]
Then according to Definition \ref{defoper1} $f\oplus g$ is the
constant zero function.

The proof of the associative and distributive laws is an easy
application of the techniques derived in the next section and will
be proved there.
\end{proof}

\section{Extension property and an alternative definition of the operations on
$\mathbb{H}(\Omega)$}

Let $D$ be a dense subset of $\Omega$. Extending a function $f$
defined on $D$ to $\Omega$ while preserving its properties (e.g.
linearity, continuity) is an important issue in functional
analysis. Recall that if $f$ is continuous on $D$ it does not
necessarily have a continuous extension on $\Omega$. Hence the
significance of the next theorem which was proved in
\cite{Sozopol}.
\begin{theorem}\label{textension}
Let $\varphi\in\mathbb{H}(D)$, where $D$ is a dense subset of
$\Omega$. Then there exists a unique $f\in\mathbb{H}(\Omega)$,
such that $f(x)=\varphi(x)$, $x\in D$. Namely,
$f=F(D,\Omega,\varphi)$.
\end{theorem}

For every two functions $f,g\in\mathbb{H}(\Omega)$ denote
$D_{fg}=\{x\in\Omega:w(f(x))=w(g(x))=0\}$. As shown already the
point-wise sum and product of H-continuous functions is not always
H-continuous. However, the restrictions of $f$, $g$, $f+g$ and
$f\times g$ on the set $D_{fg}$ are all real continuous functions,
see Theorem \ref{tDfdense}. As usual these restrictions are
denoted by $f|_{D_{fg}}$, $(f+g)|_{D_{fg}}$, etc. Using that the
set $D_{fg}$ is dense in $\Omega$ the following definition of the
operations $\oplus$ and $\otimes$ is suggested.
\begin{definition} \label{defoper2} Let $f,g\in\mathbb{H}(\Omega)$. Then
\begin{itemize}
\item[(a)] $f\oplus g$ is the unique H-continuous extension of
$(f+g)|_{D_{fg}}$ on $\Omega$ given by Theorem $\ref{textension}$,
that is, $f\oplus g=F(D_{fg},\Omega,f+g)$.%
\item[(b)] $f\otimes g$ is the unique H-continuous extension of
$(f\times g)|_{D_{fg}}$ on $\Omega$ given by Theorem
$\ref{textension}$, that is, $f\otimes g=F(D_{fg},\Omega,f\times
g)$.
\end{itemize}
\end{definition}

In order to justify the use of the notations $\oplus$ and
$\otimes$ let us immediately prove that Definitions \ref{defoper1}
and \ref{defoper2} are equivalent. Indeed, using the property
(\ref{Finclmon3}) and the fact that $f+g$ is S-continuous, see
Theorem \ref{tf*gScont}, we have
\[F(D_{fg},\Omega,f+g)(x)\subseteq F(f+g)(x)=(f+g)(x), x\in\Omega.
\]
Hence $F(D_{fg},\Omega,f+g)$ is the unique Hausdorff continuous
function satisfying the inclusion required in Definition
\ref{defoper1} which implies that $F(D_{fg},\Omega,f+g)$ is the
sum of $f$ and $g$ according to Definition \ref{defoper1}. In a
similar way one can show that $F(D_{fg},\Omega,f\times g)$ is the
product of $f$ and $g$ according to Definition \ref{defoper1}
Therefore Definitions \ref{defoper1} and \ref{defoper2} are
equivalent.

Definition \ref{defoper2} is particularly useful for evaluating
arithmetical expressions involving more than two operands since
one can evaluate the expression on a set where all operands assume
point values and then extend the answer to $\Omega$. We will
explain the procedure in detail. Let $E(+,\times,z_1,z_2,...,z_k)$
be an expression involving the operations $+$ and $\times$ and $k$
operands. Let
$\mathcal{F}=\{f_1,f_2,...,f_k\}\subset\mathbb{H}(\Omega)$. Then
the functions
\begin{equation}\label{E+x}
E(+,\times,f_1,f_2,...,f_k)\ \mbox{ and }\
E(\oplus,\otimes,f_1,f_2,...,f_k)
\end{equation}
are both well defined, the first one being S-continuous, the
second one being H-continuous. According to Theorem \ref{tDfdense}
the set $D_{\mathcal{F}}$ given by (\ref{DF}) is dense in
$\Omega$. Then a simple connection between the functions
(\ref{E+x}) is given in the next theorem, which is proved easily
using the definition of the operations and the extension property.
\begin{theorem}\label{tE+x=E+x}
For any set of functions
$\mathcal{F}=\{f_1,f_2,...,f_k\}\subset\mathbb{H}(\Omega)$ and
arithmetical expression $E(+,\times,z_1,z_2,...,z_k)$ we have
\[
F(D_{\mathcal{F}},\Omega,E(+,\times,f_1,f_2,...,f_k))=
E(\oplus,\otimes,f_1,f_2,...,f_k)
\]
\end{theorem}

As an application of Theorem \ref{tE+x=E+x} we will show that the
associative and distributive laws of a ring hold true on
$\mathbb{H}(\Omega)$ with respect to the operations $\oplus$ and
$\otimes$. This will complete the proof of Theorem \ref{tring}.

{\bf Proof of the associative and distributive laws on
$\mathbb{H}(\Omega)$.} Let $f,g,h\in\mathbb{H}(\Omega)$ and let
$D_\mathcal{F}$ be the dense subset of $\Omega$ given by
(\ref{DF}) with $\mathcal{F}=\{f,g,h\}$. Since the values of $f$,
$g$ and $h$ on the set $D_\mathcal{F}$ are all real numbers (point
intervals) we have the functions $f$, $g$ and $h$ satisfy the
associative and distributive laws with respect to the operations
$+$ and $\times$. Then using Theorem \ref{tE+x=E+x} we have
\begin{eqnarray*}
&&\hspace{-5mm}(f\oplus g)\oplus h=F(D_\mathcal{F},\Omega,(f+g)+h)=F(D_\mathcal{F},\Omega,f+(g+h))=f\oplus (g\oplus h),\\
&&\hspace{-5mm}(f\otimes g)\oplus
h=F(D_\mathcal{F},\Omega,(f\times g)\times
h)=F(D_\mathcal{F},\Omega,f\times (g\times h))
=f\otimes (g\otimes h),\\
&&\hspace{-5mm}(f(x)\oplus g(x))\otimes
h(x)=F(D_\mathcal{F},\Omega,(f+g)\times
h)\\
&&\hspace{3.4 cm}=F(D_\mathcal{F},\Omega,f\times h+g\times
h)=f\otimes h+g\otimes h,
\end{eqnarray*}
which shows that both associative laws as well as the distributive
law hold true.

\section{Definition of the operations on $\mathbb{H}(\Omega)$ through the order convergence structure}

Partial order can be defined for intervals in different ways. Here
we consider the partial order on $\mathbb{IR}$ given by
\begin{equation}\label{iorder}
[\underline{a},\overline{a}]\leq [\underline{b},\overline{b}]\
\Longleftrightarrow\ \underline{a}\leq \underline{b},\
\overline{a}\leq \overline{b}\ .
\end{equation}
The partial order on $\mathbb{H}(\Omega)$ which is induced by
(\ref{iorder}) in a  point-wise way, that is, for
$f,g\in\mathbb{H}(\Omega)$
\begin{equation}\label{forder}
f\leq g\ \Longleftrightarrow \ f(x)\leq g(x),\ x\in\Omega,
\end{equation}
is naturally associated with $\mathbb{H}(\Omega)$. For example, it
was shown in \cite{Anguelov} that the set $\mathbb{H}(\Omega)$ is
Dedekind order complete with respect to this order.

The order convergence of sequences on a poset is defined through
the partial order.

\begin{definition}\label{deforderconv}
Let $P$ be a poset with a partial order $\leq$. A sequence
$(f_n)_{n\in\mathbb{N}}$ on $P$ is said to order converge to $f\in
P$ if there exist on $P$ an increasing sequence
$(\alpha_n)_{n\in\mathbb{N}}$ and a decreasing sequence
$(\beta_n)_{n\in\mathbb{N}}$ such that $\alpha_n\leq \xi_n\leq
\beta_n$, $n\in\mathbb{N}$, and
$f=\sup_{n\in\mathbb{N}}\alpha_n=\inf_{n\in\mathbb{N}}\beta_n$.
\end{definition}

It is well known that in general the order convergence on a poset
is not topological in the sense there there is no topology with
class of convergent sequences exactly equal to the class of the
order convergent sequences. In particular this is the case of
$\mathbb{H}(\Omega)$ with respect to the partial order
(\ref{forder}). However, the order convergence induces on
$\mathbb{H}(\Omega)$ the structure of the so called FS sequential
convergence space. See \cite{Beattie} for the the definition and
further details on FS sequential convergence spaces and
convergence (filter) spaces.

The concept of Cauchy sequence in general can not be defined
within the realm of sequential convergence only but rather using
the stronger concept of a convergence space, \cite{Beattie}. It
was shown in \cite{Oconv} that the order convergence structure on
$C(\Omega)$) is induced by a convergence space and we have the
following characterization of the Cauchy sequences. Let
$(f_n)_{n\in\mathbb{N}}$ be a sequence on $C(\Omega)$. Then
\begin{equation}\label{Cauchychar}
(f_n)_{n\in\mathbb{N}}\ \mbox{ is Cauchy}\ \Longleftrightarrow\
\left\{\begin{tabular}{l}There exists a decreasing sequence\\
$(\beta_n)_{n\in\mathbb{N}}$ on $C(\Omega)$ such that
$\inf\beta_n=0$\\ and $f_m-f_k\leq\beta_n$, $m,k\geq n$,
$n\in\mathbb{N}$.\end{tabular}\right.
\end{equation}
It was also shown in \cite{Oconv} that the convergence space
completion of $C(\Omega)$ is $\mathbb{H}(\Omega)$. More precisely,
we have the following theorem.
\begin{theorem}\label{tcompletion}
\begin{itemize} \item[(i)] For every Cauchy sequence
$(f_n)_{n\in\mathbb{N}}$ on $C(\Omega)$ there exists
$f\in\mathbb{H}(\Omega)$ such that $f_n\rightarrow f$. \item[(ii)]
For every $f\in\mathbb{H}(\Omega)$ there exists a Cauchy sequence
$(f_n)_{n\in\mathbb{N}}$ on $C(\Omega)$ such that $f_n\rightarrow
f$. Moreover, the sequence $(f_n)_{n\in\mathbb{N}}$ can be
selected to be either increasing or decreasing.
\end{itemize}
\end{theorem}

\begin{definition}\label{defoper3}
Let $f,g\in\mathbb{H}(\Omega)$ and let $(f_n)_{n\in\mathbb{N}}$
and $(g_n)_{n\in\mathbb{N}}$ be the Cauchy sequences on
$C(\Omega)$ existing in terms of Theorem \ref{tcompletion}, that
is, we have $f_n\rightarrow f$, $g_n\rightarrow g$. Then $f\oplus
g$ is the order limit of $(f_n+g_n)_{n\in\mathbb{N}}$ and
$f\otimes g$ is the order limit of $(f_n\times
g_n)_{n\in\mathbb{N}}$.
\end{definition}
Let us first note that the order limits stated in Definition
\ref{defoper3} do exist. Indeed, since $(f_n)_{n\in\mathbb{N}}$
and $(g_n)_{n\in\mathbb{N}}$ are Cauchy sequences on $C(\Omega)$
one can see from (\ref{Cauchychar}) that their sum
$(f_n+g_n)_{n\in\mathbb{N}}$ and their product $(f_n\times
g_n)_{n\in\mathbb{N}}$ are Cauchy sequences. Hence according to
Theorem \ref{tcompletion}(i) they both order converge.

To establish the consistency of the Definition \ref{defoper3} we
need to show that $f\oplus g$ and $f\otimes g$ do not depend on
the particular choice of the sequences. Let
$(f_n^{(1)})_{n\in\mathbb{N}}$, $(f_n^{(2)})_{n\in\mathbb{N}}$,
$(g_n^{(1)})_{n\in\mathbb{N}}$, $(g_n^{(2)})_{n\in\mathbb{N}}$ be
Cauchy sequences on $C(\Omega)$ such that $f_n^{(i)}\rightarrow
f$, $g_n^{(i)}\rightarrow g$, $i=1,2$. We will show that the
sequences $(f_n^{(1)}+g_n^{(1)})_{n\in\mathbb{N}}$ and
$(f_n^{(2)}+g_n^{(2)})_{n\in\mathbb{N}}$ converge to the same
limit and that the sequences $(f_n^{(1)}\times
g_n^{(1)})_{n\in\mathbb{N}}$ and $(f_n^{(2)}\times
g_n^{(2)})_{n\in\mathbb{N}}$ converge to the same limit. Denote by
$(f_n)_{n\in\mathbb{N}}$ and $(g_n)_{n\in\mathbb{N}}$ the trivial
mixtures of the sequences $(f_n^{(1)})_{n\in\mathbb{N}}$,
$(f_n^{(2)})_{n\in\mathbb{N}}$ and, respectively,
$(g_n^{(1)})_{n\in\mathbb{N}}$, $(g_n^{(2)})_{n\in\mathbb{N}}$,
that is, $f_{2n-1}=f_n^{(1)}$, $f_{2n}=f_n^{(2)}$,
$g_{2n-1}=g_n^{(1)}$, $g_{2n}=g_n^{(2)}$. In an FS sequential
convergence space the trivial mixture of sequences converging to
the same limit also converges to that limit, \cite{Beattie}. Hence
we have $f_n\rightarrow f$, $g_n\rightarrow g$. It easy to see
that the sequences $(f_n)_{n\in\mathbb{N}}$,
$(g_n)_{n\in\mathbb{N}}$ are Cauchy. Hence the sequence
$(f_n+g_n)_{n\in\mathbb{N}}$ is also Cauchy, which implies that
$(f_n+g_n)_{n\in\mathbb{N}}$ converges on $\mathbb{H}(\Omega)$,
see Theorem \ref{tcompletion}(i). Therefore,
$(f_n^{(1)}+g_n^{(1)})_{n\in\mathbb{N}}$ and
$(f_n^{(2)}+g_n^{(2)})_{n\in\mathbb{N}}$, being subsequences of
the order convergent sequence $(f_n+g_n)_{n\in\mathbb{N}}$ order
converge to the same limit. In a similar way we prove that
$(f_n^{(1)}\times g_n^{(1)})_{n\in\mathbb{N}}$ and
$(f_n^{(2)}\times g_n^{(2)})_{n\in\mathbb{N}}$ converge to the
same limit.

\begin{theorem}\label{tequiv13}
Definition \ref{defoper1} and Definition \ref{defoper3} are
equivalent.
\end{theorem}
\begin{proof} Let $f,g\in\mathbb{H}(\Omega)$ and let $h_1$ be their
sum according to Definition \ref{defoper1} while $h_2$ is their
sum according to Definition \ref{defoper3}. We will show that
$h_1=h_2$. It follows from Theorem \ref{tcompletion}(ii) that we
can select increasing sequences $(f_n)_{n\in\mathbb{N}}$,
$(g_n)_{n\in\mathbb{N}}$ on $C(\Omega)$ such that $f_n\rightarrow
f$, $g_n\rightarrow g$ or equivalently,
$f=\sup\limits_{n\in\mathbb{N}}f_n$,
$g=\sup\limits_{n\in\mathbb{N}}g_n$. Then, according to Definition
\ref{defoper3}, $h_2$ is the order limit of the increasing
sequence $(f_n+g_n)_{n\in\mathbb{N}}$, that is,
$h_2=\sup\limits_{n\in\mathbb{N}}(f_n+g_n)$. On the other hand
\[
f_n(x)+g_n(x)\leq \underline{f}(x)+\underline{g}(x)\leq h_1(x),\
x\in\Omega.
\]
Therefore, $h_2\leq h_1$. In similar way by using decreasing
sequences we prove that $h_2\geq h_1$. Hence $h_1=h_2$.

The proof of the equivalence of definitions of multiplication is
done using a similar approach but is technically more complicated
an will be omitted.
\end{proof}

\end{document}